\documentclass[preprint]{imsart}

\usepackage{latexsym, amssymb, amscd, amsthm, amsxtra, amsmath,amsthm,mathtools }
\RequirePackage{natbib}
\RequirePackage[colorlinks,citecolor=blue,urlcolor=blue]{hyperref}
\usepackage{graphics, graphicx, color}
\usepackage{ifpdf}
\usepackage{multirow}
\usepackage{mathabx}
\graphicspath{{./Graphics/}{./Graphics/results/}}

\startlocaldefs

\newtheorem{thm}{Theorem}
\newtheorem{theorem}{Theorem}
\newtheorem{corollary}[thm]{Corollary}
\newtheorem{lemma}[thm]{Lemma}

\theoremstyle{definition}

\newtheorem*{defn*}{Definition}
\theoremstyle{remark}
\newtheorem{remark}{Remark}

\def\T{^{ \mathrm{\scriptscriptstyle T} }}

\newcommand{\norm}[1]{\left\Vert#1\right\Vert}

\def\Cov{\mbox{Cov}}

\def\Var{\mbox{Var}}

\def\E{\mbox{E}}

\def\half{\frac{1}{2}}

\newcommand{\epsilonv}{\mbox{\boldmath{$\epsilon$}}}

\newcommand{\Wc}{\mathcal{W}}
\newcommand{\Xc}{\mathcal{X}}

\def\1v{\mathbf 1}
\def\0v{\mathbf 0}
\def\Id{\mathbb I}

\newcommand{\bi}{\begin{itemize}}
\newcommand{\ei}{\end{itemize}}
\newcommand{\be}{\begin{enumerate}}
\newcommand{\ee}{\end{enumerate}}
\newcommand{\bb}{\begin{block}}
\newcommand{\eb}{\end{block}}
\newcommand{\ba}{\begin{align}}
\newcommand{\ea}{\begin{align}}
\newcommand{\bd}{\begin{align*}}
\newcommand{\ed}{\begin{align*}}

\def\varphi{\lambda}

\endlocaldefs

\begin{document}

\begin{frontmatter}

\title{Adjusting systematic bias in high dimensional principal component scores}
\runtitle{Bias in principal component scores}

\author{\fnms{Sungkyu} \snm{Jung}\ead[label=e2]{sungkyu@pitt.edu}}
\address{Department of Statistics, Seoul National University\\ \printead{e2}}
 \runauthor{S. Jung}

\begin{abstract}
Principal component analysis continues to be a powerful tool in dimension reduction of high dimensional data. We assume a variance-diverging model and use the high-dimension, low-sample-size asymptotics to show that even though the principal component directions are not consistent, the sample and prediction principal component scores can be useful in revealing the population structure. We further show that these scores are biased, and the bias is asymptotically  decomposed into rotation and scaling parts. We propose methods of bias-adjustment that are shown to be consistent and work well in the high dimensional situations with small sample sizes. The potential advantage of bias-adjustment is demonstrated in a classification setting.
\end{abstract}





\end{frontmatter}

\section{Introduction}

Principal component analysis (PCA) is a workhorse method of multivariate analysis, and has been used in a variety of fields for dimension reduction, visualization and as exploratory analysis. The standard estimates of principal components, obtained by either the eigendecomposition of the sample covariance matrix or the singular value decomposition of the data matrix, are now well-known to be inconsistent when the number of variables, or the dimension $d$, is much larger than the sample size $n$ \citep{Paul2007,Johnstone2009,Jung2009a}. These observations were paralleled with a vast amount of proposals on, e.g., sparse principal component estimations \citep[\emph{cf}. most notably, ][]{Zou2006}, which perform better in some models with high dimensions. 

However, the standard estimates of principal components (PCs) continue to be useful, partly due to fast computations available \citep[see, e.g.,][]{abraham2014fast}. Many of the sparse estimation methods, unfortunately, do not computationally scale well for large data with hundreds of thousands of variables. Moreover, the standard estimation has shown to be useful in some application areas such as imaging, genomics and big-data analysis \citep{doi:10.1093/nsr/nwt032}. In these areas, the sample and prediction PC scores (the projection scores of the data points onto the PC directions) are often used in the next stage of analysis. 


The prediction of PC scores has considerable  practical utility in modern data analysis.
A prominent example where the ``sample" and ``prediction" PC scores are used is the \emph{PC regression}. In particular, for
prediction and cross-validation for PC regression, the PC scores are used as explanatory variables. For prediction of the response from a new set of observations, the predicted PC scores are needed \citep{jackson2005user}. As an example, \cite{qiuhua2014use} used a PC regression in prediction of an phytoplankton abundance index. In the same vein, \emph{classification} rules are often estimated for  dimension-reduced data sets. As an instance, in forensic science, residue features from various black ballpoint inks are dimension-reduced (via PCA) then classified, based on a lab data set. New features from the field are classified using their prediction scores as an input for the classification rule \citep{adam2008classification}. 
As a more involved example, ancestry estimation in genetic association studies uses the sample PC scores obtained from a reference genotyped sample, often from large-scale public sequencing data sets \citep{zhan2013identification,
    marcus2020genetic,wang2015improved}.
    The prediction PC scores of a new sample is then matched  to the sample PC scores, in order to infer the new samples's ancestry membership \citep{zhang2020fast}.


In this paper, we revisit the standard estimates of principal components in ultra-high dimensions and reveal that while the component {directions} and {variances} are inconsistent, the sample and prediction scores are useful for moderately large sample size. For low sample sizes, the scores are biased. We quantify the bias, decompose it into two systematic parts, and propose to estimate bias-adjustment factors.

\begin{figure}[tb!]
  \centering
  \vskip -4cm
  \includegraphics[width=1\textwidth]{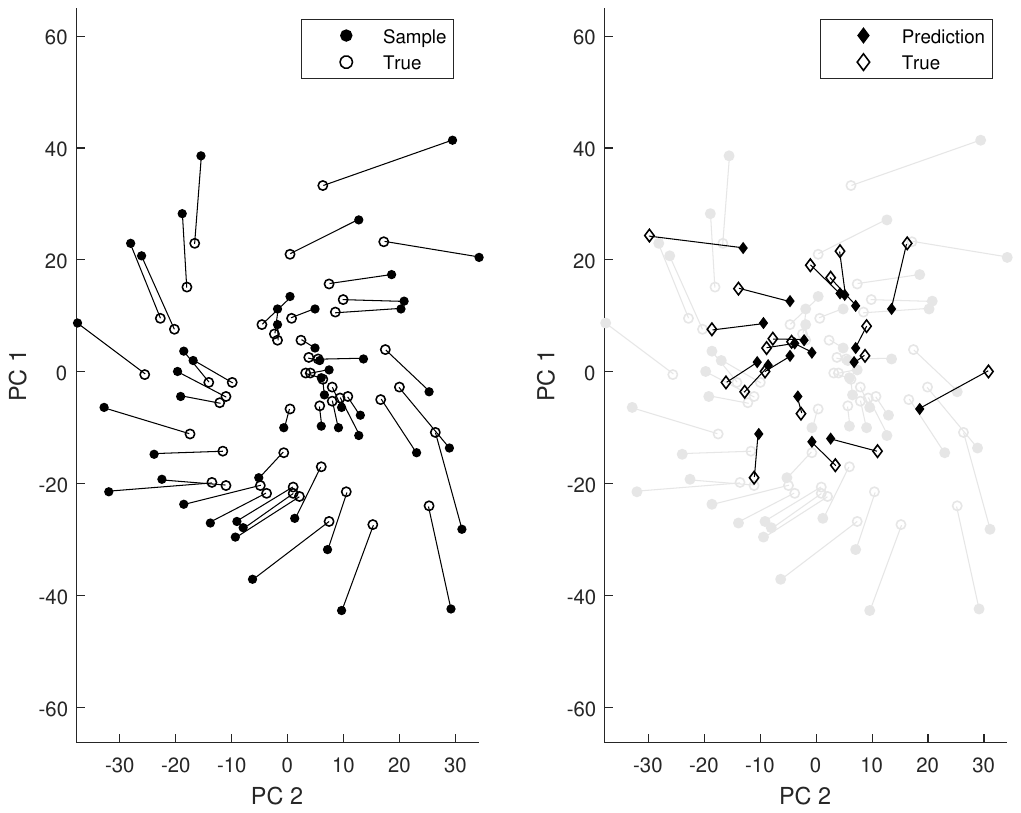}
  \vskip -4cm
  \caption{Sample and prediction principal component scores connected to their true values. This toy data set of size $(d,n) = (10000, 50)$ is generated from the spike model with $m=2$ spikes, with polynomially-decreasing eigenvalues with $\beta = 0.3$; see Section~\ref{sec:SIMmodels} for details.
  \label{fig:projectionScores2}}
\end{figure}

As a visual example of the systematic bias, a toy data set with 2 distinguishable principal components is simulated and plotted in Fig.~\ref{fig:projectionScores2}. Each observation in the data set consists of $d = 10,000$ variables. The first two sample principal component directions are estimated from $n = 50$ observations, and are used to obtain the sample and prediction scores (the latter are computed from 20 new observations). The true principal scores  are also plotted and connected to their empirical counterparts. This example visually reveals that the sample scores are systematically biased, that is, \emph{uniformly rotated} and \emph{stretched}. What is more surprising is that the prediction scores are also uniformly rotated, by  the same angle as the sample scores, and uniformly shrunk.

\begin{figure}[tb!]
  \centering
  \vskip -4cm
  \includegraphics[width=1\textwidth]{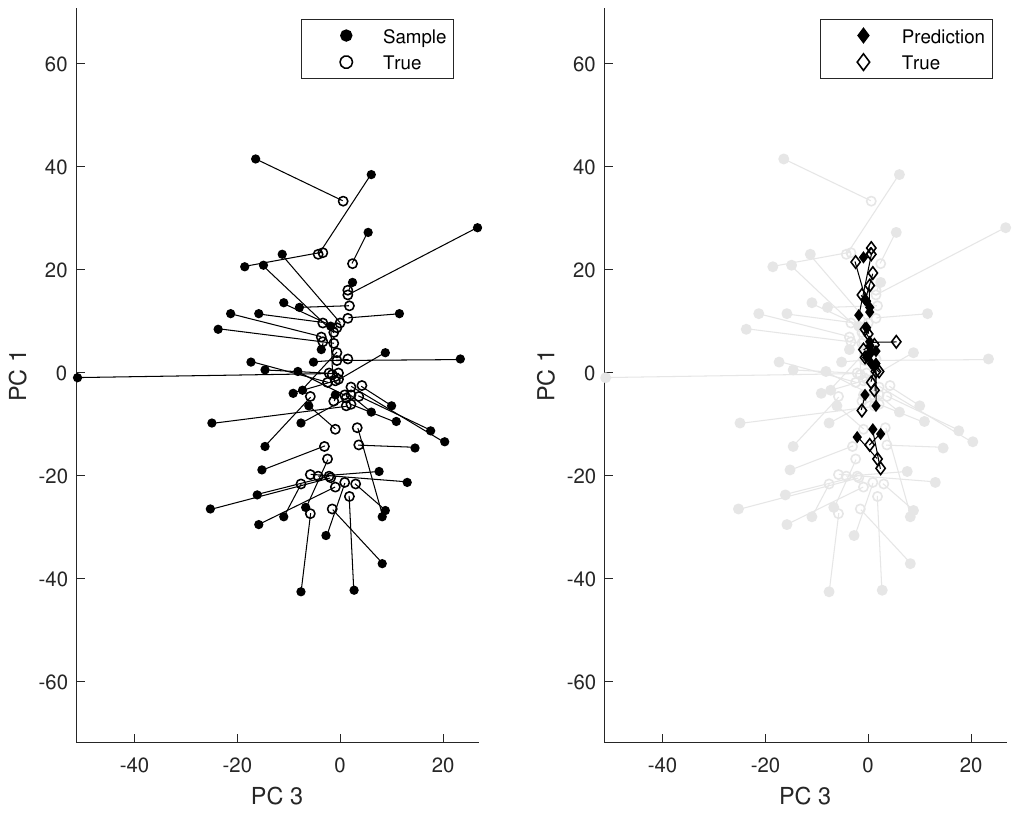}
  \vskip -4cm
  \caption{Sample and prediction principal component scores connected to their true values. Models and data are the same as in Fig.~\ref{fig:projectionScores2}.
  \label{fig:projectionScores5}}
\end{figure}

On the other hand, the third component scores from this example appear to be quite arbitrary; see Fig.~\ref{fig:projectionScores5}. (The estimates for component 3 in this example is only as good as random guess.) Moreover, unlike the first two components plotted in Fig.~\ref{fig:projectionScores2}, the sample scores of the third component are grossly inflated, while the prediction scores are much smaller than the sample scores.

In Section~\ref{sec:scores}, we provide theoretical justification of the phenomenon observed in Figs.~\ref{fig:projectionScores2} and \ref{fig:projectionScores5}, and asymptotically quantify the two parts of the systematic bias. We assume $m$-component models with diverging variances, and use the high-dimension, low-sample-size asymptotic scenario (i.e. $d\to\infty$ while $n$ is fixed). These models and asymptotics are used in giving the contrasting results of the sample and prediction scores. The correlation coefficients between the sample (or prediction) and true scores turn out to be close to 1, for large signals and large sample sizes, indicating the situations where the principal component scores are most useful.

Since the bias is asymptotically quantified, the natural next step is to adjust the bias by estimating the bias-adjustment factor. In Section~\ref{sec:bias-adjustment}, we propose a simple, yet consistent, estimator and several variants of estimators based on the idea of Jackknife. Adjusting these biases improves the performance of prediction modeling, and we demonstrate its potential by an example involving classification. Results from numerical studies are summarized in Section~\ref{sec:numericalstudies}.


There are several related works on the principal component scores in high dimensions \citep{lee2010convergence,fan2013large,lee2014convergence,sundberg2016exploratory,shen2016statistics,SJOS:SJOS12264,wang2017asymptotics,Jung2017}. This paper is built upon these previous findings.
In particular, this paper is a continuation of the author's previous work  \citep{Jung2017}, and intermediate results are borrowed from there. While the scaling and rotation of the sample scores were previously identified in \cite{Jung2017} as well as in \cite{SJOS:SJOS12264}, the main contributions of this paper are  \emph{i)} the quantification of the asymptotic bias for the \emph{prediction} scores, which has not been addressed, and \emph{ii)} a consistent estimation of the bias-adjustment factor. Under the ``random-matrix'' asymptotic scenario, i.e., $d/n \to c \in (0, \infty)$, \cite{lee2010convergence} discussed a bias adjustment of principal component scores. Our work extends \cite{lee2010convergence} to the high-dimension, low-sample-size asymptotic scenario. Note that the asymptotic \emph{rotational} bias was not identified in \cite{lee2010convergence}, due to larger sample size $n \asymp d$ considered there.
A survey of high-dimension, low-sample-size asymptotics can be found in \cite{AOSHIMA2017}.


\section{Asymptotic behavior of principal component scores}\label{sec:scores}

\subsection{Model and assumptions}

Let $\Xc = [X_1,\ldots,X_n]$ be a $d \times n$ data matrix, where each $X_i$ is mutually independent and has zero mean and covariance matrix $\Sigma_d$. Population principal components are obtained by the eigendecomposition of $\Sigma_d = U\Lambda U\T$, where $\Lambda = \mbox{diag}(\lambda_{1},\ldots, \lambda_{d})$ is the diagonal matrix of principal component variances and $U = [u_1,\ldots, u_d]$ consists of principal component directions. For a fixed $m$, we assume an $m$-component model, where the first $m$ component variances are distinguishably larger than the rest. Specifically, the larger variances increase at the same rate as the dimension $d$, i.e. $\lambda_i \asymp d$, which was previously noted as the ``boundary situation'' \citep{Jung2012}. This diverging-variance condition seems to be more realistic than the other simpler cases $\lambda_i \gg d$ (i.e., $\lambda_i / d \to \infty$) and $\lambda_i \ll d$ \citep{SJOS:SJOS12264, shen2016statistics}, and is satisfied for high-dimensional models used in factor analysis \citep{fan2013large,li2016embracing,sundberg2016exploratory}. In a more general asymptotic scenario of $d/n \to \infty$, our condition, $\lambda_i \asymp d$, is akin to the condition, $\lim_{n\to\infty}d/(n\lambda_i) = c_i \in (0,\infty)$, assumed in \cite{shen2016statistics} and \cite{wang2017asymptotics}. In particular, in the \emph{ultra-high dimensional} case of $n \asymp \log(d)$, as defined in \cite{fan2008sure}, we have $d^{1-\epsilon} \ll d/n \ll d^{1+\epsilon}$ for any $\epsilon >0$. Thus, although not identical, the assumption $ \lambda_i\asymp d/n $ of \cite{shen2016statistics} and \cite{wang2017asymptotics} is similar to (A1) below, $\lambda_i\asymp d$, in the ultra-high dimensional case.

We assume that the population principal component variances satisfy the following:
\begin{itemize}
\item[(A1)] $ \lambda_i = \sigma_i^2 d, \ i = 1,\ldots,m, \ \sigma_1^2 \ge \cdots \ge \sigma_m^2$.
\item[(A2)] $\lim_{d\to\infty} \sum_{i=m+1}^d \lambda_i /d := \tau^2 \in (0, \infty)$.
\item[(A3)] 
     There exists $B <\infty$ such that for all $i>m$, $\limsup_{d\to\infty} \lambda_i < B$.
\end{itemize}

The conditions (A2) and (A3) are used to allow $\lambda_i$ for $i>m$ increase as $d$ increases. All of our results hold when the condition (A3) is relaxed to, e.g., allow the situation that $\lambda_i \asymp d^{\alpha}$, $\alpha < 1/2$. Such generalization is straightforward, but invites nonintuitive technicality \citep[see, e.g., ][]{Jung2012,Jung2017}.
By decomposing each independent observation into the first $m$ components and the remaining term, we write
\begin{equation}\label{eq:X_j-decomposition}
X_j = \sum_{i=1}^m \lambda_i^{1/2} u_i z_{ij} + \sum_{i=m+1}^d \lambda_i^{1/2} u_i z_{ij},\quad (j = 1,\ldots,n),
\end{equation}
where $z_{ij}$ is the normalized principal component score.

\bi
 \item[(A4)] For each $j = 1,2,\ldots$, $(z_{1j},z_{2j},\ldots)$ is a sequence of independent random variables such that for any $i$, $\E(z_{ij}) = 0$, $\Var(z_{ij}) = 1$, and that the fourth moment of $z_{ij}$ is uniformly bounded.
\ei

\subsection{Sample and prediction principal component scores}
Suppose we have a data matrix $\Xc = [X_1,\ldots,X_n]$ and a vector $X_*$, independently drawn from the same population with principal component directions $u_i$. The principal component analysis is performed for data $\Xc$ and is used to predict the principal component scores of $X_*$.

We define the  $i$th \emph{true principal component scores} of $\Xc$ as the vector of $n$ projection scores:
\begin{equation}\label{eq:truePCscore}
w_i\T = u_i\T \Xc = (w_{i1},\ldots,w_{in}), \quad (i=1,\ldots,d),
\end{equation}
where $w_{ij} = u_i\T X_j = \sqrt{\lambda_i}z_{ij}$. The last equality is given by the decomposition of $X_j$ in (\ref{eq:X_j-decomposition}).
Likewise, the true $i$th principal component score of $X_*$ is $w_{i*} = u_i\T X_* = \sqrt{\lambda_i} z_{i*}$.

The classical estimators of the pair of the $i$th principal component direction and variance are $(\hat{u}_i,\hat\lambda_i)$, obtained by either the eigendecomposition of the sample covariance matrix $S_d = n^{-1}\Xc\Xc^T$,
$$ S_d = \sum_{i=1}^n \hat\lambda_i \hat{u}_i \hat{u}_i\T, $$
or by the singular value decomposition of the data matrix,
\begin{equation} \label{eq:SVD on X}
\Xc = \sqrt{n} \sum_{i=1}^n  \sqrt{\hat\lambda_i} \hat{u}_i \hat{v}_i\T,
\end{equation}
where $\hat{v}_i$ is the right singular vector of $\Xc$. By replacing $u_i$ in (\ref{eq:truePCscore}) with its estimator $\hat{u}_i$, we define the $i$th \emph{sample principal component scores} of $\Xc$ as
\begin{equation}\label{eq:samplePCscore}
\hat{w}_i\T = \hat{u}_i\T \Xc = (\hat{w}_{i1},\ldots,\hat{w}_{in}), \quad (i=1,\ldots,n).
\end{equation}
The sample principal component scores are in fact weighted right singular vectors of $\Xc$; comparing to (\ref{eq:SVD on X}), $\hat{w}_i = \sqrt{n\hat\lambda_i} \hat{v}_i$.

For an independent observation $X_*$, the definition (\ref{eq:samplePCscore}) gives
$$
\hat{w}_{i*} = \hat{u}_i\T X_*,
$$
which is called the $i$th \emph{prediction principal component score} for $X_*$.

\subsection{Main results}\label{sec:mainresults}

Denote  $W_1 =  (\sigma_i z_{ij})_{i,j} = (d^{-1/2} w_{ij})_{i,j} = d^{-1/2}[u_1,\ldots,u_m]\T \Xc$
for the $m \times n$ matrix of the scaled true scores for the first $m$ principal components. The $i$th row of $W_1$ is $d^{-1/2}w_i\T$. Similarly, the scaled sample scores for the first $m$ principal components are denoted by
$\widehat{W}_1 = d^{-1/2} [\hat{u}_1,\ldots,\hat{u}_m]\T \Xc$.

For a new observation $X_*$, write $W_{*} = d^{-1/2} ({w}_{1*},\ldots,{w}_{m*})\T$
and $\widehat{W}_{*} = d^{-1/2} (\hat{w}_{1*},\ldots,\hat{w}_{m*})\T$ for the scaled true scores and prediction scores, respectively, of the first $m$ principal components.

Write $\Wc = W_1 W_1\T $  for  the scaled $ m\times m$ sample covariance matrix of the first $m$  scores.
Let $\{\varphi_i(S), v_i(S)\}$ denote the $i$th largest eigenvalue-eigenvector pair of a non-negative definite matrix $S$ and $v_{ij}(S)$ denote the $j$th loading of the vector $v_i(S)$.  For a sequence $A_d$ of random matrices, we say $A_d = O_p(b_d)$ if all elements of $A_d/b_d$ are uniformly stochastically bounded. Note that $A_d = O_p(1)$ implies $\|A_d\|_F = O_p(1)$.

\begin{theorem} \label{thm:main}
Assume the $m$-component model under Conditions (A1)--(A4) and let $n > m \ge 0$ be fixed and $d\to\infty$. Then, the first $m$ sample and prediction scores are systematically biased:
    \begin{align}
     \widehat{W}_1 &= S R\T W_1  + O_p(d^{-1/4}), \label{eq:scarot1}\\
     \widehat{W}_{*} &= S^{-1} R\T  W_*  + O_p(d^{-1/2}),  \label{eq:scarot2}
     \end{align}
    where $R = [v_1(\Wc),\ldots, v_m(\Wc)]$, $S = \mbox{diag}(\rho_1,\ldots,\rho_m),$  
and $\rho_k = \sqrt{ 1+ \tau^2/\varphi_k(\Wc) }$.
 Moreover,
 for $k > m$,
     \begin{align}
     \hat{w}_{kj} &= O_p(d^{1/2}), \quad j = 1,\ldots,n, \label{eq:score-sample-diverge}\\
      \hat{w}_{k*} &= O_p(1).  \label{eq:score-prediction-stable}
     \end{align}
\end{theorem}
Our main results show that the first $m$ sample and prediction scores are comparable to the true scores.
The asymptotic relation tells that for large $d$, the first $m$ sample scores in $\widehat{W}_1$ converge to the true scores in $W_1$, uniformly rotated and scaled for all data points.
It is thus valid to use the first $m$ sample principal scores for exploration of important data structures, and to reduce the dimension of the data space from $d$ to $m$ in the high-dimension, low-sample-size context.

Theorem~\ref{thm:main} explains and quantifies the two parts of the bias, exemplified in Fig.~\ref{fig:projectionScores2}. In particular, the same rotational bias applies to both sample and prediction scores. The scaling bias factors $\rho_k$ in the matrix $S$ are all greater than 1. Thus, while the sample scores are all stretched, the prediction scores are all shrunk.
The second part of the theorem shows that the magnitude of inflation for the  sample scores of the ``noise'' component (see, e.g., component 3 scores in Fig.~\ref{fig:projectionScores5}) is of order $d^{1/2}$. On the other hand, the prediction scores of the noise component do not diverge.
%

\begin{remark}
Suppose $m = 1$ in Theorem~\ref{thm:main}. Then the sample and prediction scores are simply proportionally-biased in the limit: $\hat{w}_{1j}/w_{1j} \to \rho_1$ and $\hat{w}_{1*}/w_{1*} \to \rho^{-1}_1$ in probability as $d\to\infty$. 
\end{remark}
\begin{remark}
Suppose that the limit $n\to\infty$ is taken for the expression (\ref{eq:scarot1})  and (\ref{eq:scarot2}). Then from the classical asymptotic results on the $m\times m$ covariance matrix $\Wc$ \citep[\emph{cf.}][]{anderson1963asymptotic}, $S = I_m + O_p(\frac{1}{n})$ and $R = I_m + O_p(\frac{1}{n})$. That is, in the limit $d\to\infty$, the limiting bias is of order $n^{-1}$.
\end{remark}

The proof of Theorem~\ref{thm:main} relies on the asymptotic behavior of the principal component direction and variance, which is now well-understood; see  \cite{Jung2017} for the asymptotic regime of $d\to\infty$, $n$ fixed; \cite{shen2016statistics} and \cite{wang2017asymptotics} for the asymptotic regime of $d\to\infty$, $n\to\infty$ and $d/n \to \infty$. For reference we restate it here.
\begin{lemma}\label{thm:PCA-asymptotic} [Theorem S2.1, \cite{Jung2017}]
 Assume the conditions of Theorem~\ref{thm:main}.
  (i)  the sample principal component variances converge in probability as $d \to \infty$;
   \begin{equation*}
  d^{-1} n \hat\lambda_{i} = \left\{
                                   \begin{array}{ll}
                                     \varphi_i (\Wc) + \tau^2 + O_p(d^{-1/2}) , & i = 1,\ldots,m; \\
                                     \tau^2 + O_p(d^{-1/2}), & i = m+1,\ldots,n.
                                   \end{array}
                                 \right.
   \end{equation*}
  (ii) The inner product between sample and population PC directions converges  in probability as $d \to \infty$;
   \begin{equation*}
   \hat{u}_{i}\T u_{j} = \left\{
                                     \begin{array}{ll}
                                       \rho_i^{-1} v_{ij}(\Wc) + O_p(d^{-1/2}) , & i,j = 1,\ldots,m; \\
                                       O_p(d^{-1/2}), & \mbox{ otherwise.}
                                     \end{array}
                                   \right.
   \end{equation*}
\end{lemma}
This result is abridged later in Section~\ref{sec:PCAasymp} for discussion. To handle prediction scores, we need in addition the following observation, summarized in Lemma~\ref{lem:error_term_bound}. For each $k = 1,\ldots,m$, the $k$th projection score $\hat{w}_{k*}$ is decomposed into
\begin{align}\label{eq:projectionNew}
\hat{w}_{k*} = \hat{u}_k\T X_{*} = \sum_{i=1}^m w_{i*} \hat{u}_k\T u_i + \epsilon_{k*},
\end{align}
where $\epsilon_{k*} =  \sum_{i=m+1}^d w_{i*} \hat{u}_k\T  u_i.$
In the next lemma, we show that the ``error term,'' $\epsilon_{k*}$,  is stochastically bounded.

\begin{lemma}\label{lem:error_term_bound}
 Assume the $m$-component model with (A1)--(A4) and let $n > m \ge 0$ be fixed. For $k = 1,\ldots,n$, $\E( \epsilon_{k*} | W_1) = 0$, and
   \begin{align} \label{eq:variancelimit}
    \lim_{d\to\infty}\Var ( \epsilon_{k*}   \mid W_1)
    &= \upsilon^2_O  / (\varphi_k(\Wc) + \tau^2 ), \quad \hbox{for $k \le m$;} \\
\lim_{d\to\infty} \frac{1}{n-m} \sum_{k = m+1}^n \Var ( \epsilon_{k*}   \mid W_1)
    &= \upsilon^2_O  / \tau^2   ,  \label{eq:variancelimit2}
\end{align}
where $\upsilon^2_O = \lim_{d\to\infty} d^{-1}\sum_{i=m+1}^d \lambda^2_i$. As $d\to\infty$,  $\epsilon_{k*}  = O_p(1)$.
\end{lemma}

Lemmas \ref{thm:PCA-asymptotic} and \ref{lem:error_term_bound} facilitate an interpretation of the results in Theorem~\ref{thm:main}.
Intuitively, the overestimation of the sample principal variances, in Lemma~\ref{thm:PCA-asymptotic}(i), causes the sample scores to be stretched, while the inconsistency of $\hat{u}_i$ leads to smaller $\hat{u}_i\T u_i$ in Lemma~\ref{thm:PCA-asymptotic}(ii), which then results in the deflation of the projection scores  (\ref{eq:projectionNew}).
Proofs of Theorem~\ref{thm:main} and all other results can be found in the supplementary material.

Next result shows that the sample and true scores (or prediction and true scores) are highly correlated with each other.
For this, we compute the inner product between the standardized sample scores $\hat{w}_k / \sqrt{\hat{w}_k\T \hat{w}_k}$ and true scores $w_k / \sqrt{{w}_k\T {w}_k}$. Define for a pair $(x,y)$ of $n$-vectors $r(x,y) = x\T y / \sqrt{x\T x \cdot y\T y}$, which is an empirical correlation coefficient between $x$ and $y$ when the mean is assumed to be zero.
%

\begin{theorem}\label{thm:correlations}
 Let $\zeta_{kj} = \lambda_k(\Wc)/(\sum_{\ell= 1}^m v^2_{\ell j}(\Wc) \lambda_\ell(\Wc) )$ and $\bar\zeta_{kj}  = \sigma^2_k/(\sum_{\ell= 1}^m v^2_{\ell j}(\Wc) \sigma^2_\ell )$.
Under the assumptions of Theorem~\ref{thm:main}, as $d\to\infty$,  for  $k,j = 1\ldots,m$,
\begin{itemize}
\item[(i)] $r(\hat{w}_k,w_j) \to v_{kj}(\Wc) \zeta_{kj}^{1/2}$ in probability ;
\item[(ii)] $ \lim_{d\to\infty}\mbox{Corr}(\hat{w}_{k*},w_{j*} \mid W_1) = v_{kj}(\Wc) \bar\zeta_{kj}{ }^{1/2}$.
\end{itemize}
\end{theorem}

\begin{remark}
In the special case, $m = 1$, both the sample and prediction scores of the first principal component are perfectly correlated with the true scores, in the limit. Specifically,  Theorem~\ref{thm:correlations} leads that
$|r(\hat{w}_1,w_1)| \to 1$ in probability and $|\mbox{Corr}(\hat{w}_{k*},w_{j*} )| \to 1$ as $d\to\infty$.
\end{remark}

\begin{remark}
The somewhat complex limiting quantity $v_{kj}(\Wc) \zeta_{kj}^{1/2}$ is an artifact of the fixed sample size. To simplify the expression for the case $k=j$, write
$$ \left( v_{kk}(\Wc) \zeta_{kk}^{1/2} \right)^2 =  \frac{1}{1+\xi_{k}(\Wc)},\quad \xi_{k}(\Wc) = \sum_{\ell \neq k} v_{\ell k}^2 (\Wc)  \frac {\lambda_\ell(\Wc) }{\lambda_k(\Wc)}.$$
Note that $\Wc = W_1 W_1\T$ is proportional to the sample covariance matrix of the first $m$ true scores, and that $v_{kk}(\Wc)$ is the inner product between the $k$th sample and theoretical principal component directions of the data set $W_1$, where the number of variables, $m$, is smaller than the sample size $n$. Therefore, we expect that $|v_{kk}(\Wc)| \approx 1$ and $\xi_k(\Wc) \approx 0$ for large sample size $n$. Taking the additional limit  $n \to\infty$, the results in Theorem~\ref{thm:correlations} become more interpretable:
$$|r(\hat{w}_k,w_j)| \to 1_{(k=j)} \mbox{ in probability, and }
|\mbox{Corr}(\hat{w}_{k*},w_{j*})| \to 1_{(k=j)},$$
as $d\to \infty, n\to\infty$ (limits are taken progressively).
\end{remark}
\begin{remark}\label{remark:conjecture}
What is the correlation coefficient $r(\hat{w}_k,w_k)$ for $k>m$ in the limit $d\to\infty$? In an attempt to answer this question, we note $\hat{w}_k = ({n\hat\lambda_k})^{1/2} \hat{v}_k$, $\hat{v}_k = v_k(\Xc\T \Xc)$ and $\Xc\T \Xc = \sum_{i=1}^d w_i w_i\T$. Thus,
$$r(\hat{w}_k,w_k) = w_k\T v_k(\sum_{i=1}^d w_i w_i\T) / \sqrt{\lambda_k},$$
and it is natural to guess that the dependence of $\hat{v}_k$ on any $w_i$, including the case $i = k$, would diminish as $d$ tends to infinity.
In fact, $d^{-1}\Xc\T \Xc$ converges to the rank-$m$ matrix $S_0 := W_1\T W_1 + \tau^2 I_n$ \citep{Jung2012}, and $w_k$ and $S_0$ are independent. Thus, it is reasonable to conjecture that $\lim_{d\to\infty}\E [r(\hat{w}_k,w_k)] = 0$, for $k>m$. Unfortunately,
in the limit $d\to\infty$, the $k$th, $k>m$, eigenvector of $d^{-1}\Xc\T \Xc$  becomes an arbitrary choice in the left null space of $W_1$. Due to this non-unique eigenvector, the inner product $w_k\T v_k(S_0)$ is not defined, and consequently discussing the convergence of $r(\hat{w}_k,w_k)$ is somewhat demanding. We numerically confirm the conjecture in Section~\ref{sec:numeric-theory}.
\end{remark}

\subsection{Inconsistency of the direction and variance estimators}\label{sec:PCAasymp}

The findings in the previous subsection may be summarized as that the first $m$ principal component scores convey about the same visual information as the true values when displayed. (The information is further honed by the bias adjustment in Section~\ref{sec:bias-adjustment}.) In a practical point of view, the scores and their graph matter the most.

On the other hand, a quite different conclusion about the standard principal component analysis is made when the estimator $\hat{u}_i$ is of interest. The asymptotic behavior of the direction $\hat{u}_i$ as well as the variance estimator $\hat\lambda_{i}$ are obtained as a special case of Lemma~\ref{thm:PCA-asymptotic}. Under our model,
\begin{equation}\label{eq:PCAinconsistency}
   (\hat{u}_{i}\T u_{i}  , d^{-1} n \hat\lambda_{i}) \to \left\{
                                   \begin{array}{ll}
                                     (\rho_i^{-1} v_{ii}(\Wc), \varphi_i (\Wc) + \tau^2 ) , & i = 1,\ldots,m; \\
                                     (0, \tau^2 ), & i = m+1,\ldots,n.
                                   \end{array}
                                 \right.
\end{equation}
in probability as $d\to\infty$ ($n$ is fixed).

The variance estimator $\hat\lambda_i$, for $ i \le m$, is asymptotically proportionally-biased.
Specifically, ${\hat\lambda_i}/{\lambda_i} \to  ({\varphi_i(\Wc)+ \tau^2})/({n\sigma_i^2})$
in probability as $d\to\infty$.
Thus by using a classical result on the expansion of the eigenvalues of $\Wc$ for large $n$,
$$E({\hat\lambda_i}/{\lambda_i}) \to 1 + \frac{1}{n}\left[\sum_{j \neq i}^{m} \frac{\sigma_j^2}{\sigma_{i}^2 - \sigma_j^2} + \frac{\tau^2}{\sigma_i^2} \right]+ O(n^{-2}),$$
as $d\to\infty$. Note that even when $m = 1$, the bias is still of order $n^{-1}$. This proportional bias may be empirically adjusted, using good estimates of $\sigma_i^2$ and $\tau^2$. We do not pursue it here. Note that all empirical principal component variances, for $i>m$, converge to $\tau^2/ n$, when scaled by $d$, and thus do not reflect any information of the population.

The result (\ref{eq:PCAinconsistency}) also shows that the direction estimator $\hat{u}_i$ is inconsistent and asymptotically-biased, compared to $u_i$. The estimator $\hat{u}_i$ is closer to $u_i$ when $\rho_i^{-1} |v_{ii}(\Wc)|$ is closer to 1. It is impossible to achieve $\rho_i^{-1} |v_{ii}(\Wc)| \to 1$ since for finite $n$, both $|v_{ii}(\Wc)|$ and $\rho_i^{-1}$ are strictly less than 1.
Although the ``angle'' between $\hat{u}_i$ and $u_i$ is quantified in (\ref{eq:PCAinconsistency}), the theorem itself is useless in adjusting the bias. This is because that the direction to which $\hat{u}_i$ moves away from $u_i$ is random, i.e. uniformly distributed; see \cite{wang2017asymptotics} for the limiting distribution of $\hat{u}_{i}$ under a general asymptotic scenario of $d/n \to \infty$, while $d/(n\lambda_i)^{-1}$ is bounded.


In short, while the bias in the principal component direction is challenging to remove, the bias in the sample and prediction scores can be  quantified and removed.



\section{Bias-adjusted scores}\label{sec:bias-adjustment}

In this section, we describe and compare several choices for the estimation of the \emph{bias-adjustment factor} $\rho_i$.
{Note that both sample and prediction scores are rotated by the same direction and amount, specified in the matrix $R$. For applications requiring score matching (e.g., classification rules trained on the sample scores or the ancestry estimation discussed in the introduction), coordinate-free methods are often used and there is less practical advantage in estimating $R$. We focus on adjusting the scores by estimating $\rho_i$.}

Suppose that the number of effective principal components, $m$, is prespecified or estimated in advance. Our first estimator is obtained by replacing $\tau^2$ and $\lambda_i(\Wc)$ in $\rho_i = \sqrt{1 + \tau^2 / \lambda_i(\Wc)}$ with reasonable estimators. In particular,
we set
\begin{equation}
\tilde\tau^2 = \frac{ \sum_{i = m+1}^n \hat\lambda_i }{ n-m} \frac{n}{d},\quad
\tilde\lambda_i(\Wc) = d^{-1} n \hat\lambda_{i} - \tau^2,
\end{equation}
and
\begin{equation}\label{eq:rho-est1}
\tilde{\rho}_i = \sqrt{1 + \tilde\tau^2 / \tilde\lambda_i(\Wc)}, \quad (i=1,\ldots,m).
\end{equation}
This simple estimator $\tilde{\rho}_i$ is in fact consistent.

\begin{corollary}\label{cor:rho_estimate} Suppose the assumptions of Lemma~\ref{thm:PCA-asymptotic} are satisfied. Let $d\to\infty$. For $i = 1,\ldots,m$, conditional to $W_1$, $\tilde\tau^2$, $\tilde\lambda_i(\Wc)$ and $\tilde{\rho}_i$ are consistent estimators of $\tau^2$, $\lambda_i(\Wc)$ and $\rho_i$, respectively.
\end{corollary}%

Using (\ref{eq:rho-est1}), the bias-adjusted sample and prediction scores are
$\hat{w}_i^{\rm (adj)} = \tilde{\rho}_i^{-1} \hat{w}_i$ and
$\hat{w}_{i*}^{\rm (adj)} = \tilde{\rho}_i \hat{w}_{i*}$ for $i = 1,\ldots,m$. The sample and prediction scores matrices in (\ref{eq:scarot1}) and (\ref{eq:scarot2}) are then adjusted to, using $\tilde{S} = \mbox{diag}(\tilde{\rho}_`,\ldots,\tilde{\rho}_m)$,
\begin{equation}\label{eq:bias-adjustment}
\widehat{W}_1^{\rm (adj)} = \tilde{S}^{-1} \widehat{W}_1, \quad \widehat{W}_{*}^{\rm (adj)} = \tilde{S} \widehat{W}_{*}.
\end{equation}

\begin{figure}[tb!]
  \centering
  \vskip -4cm
  \includegraphics[width=1\textwidth]{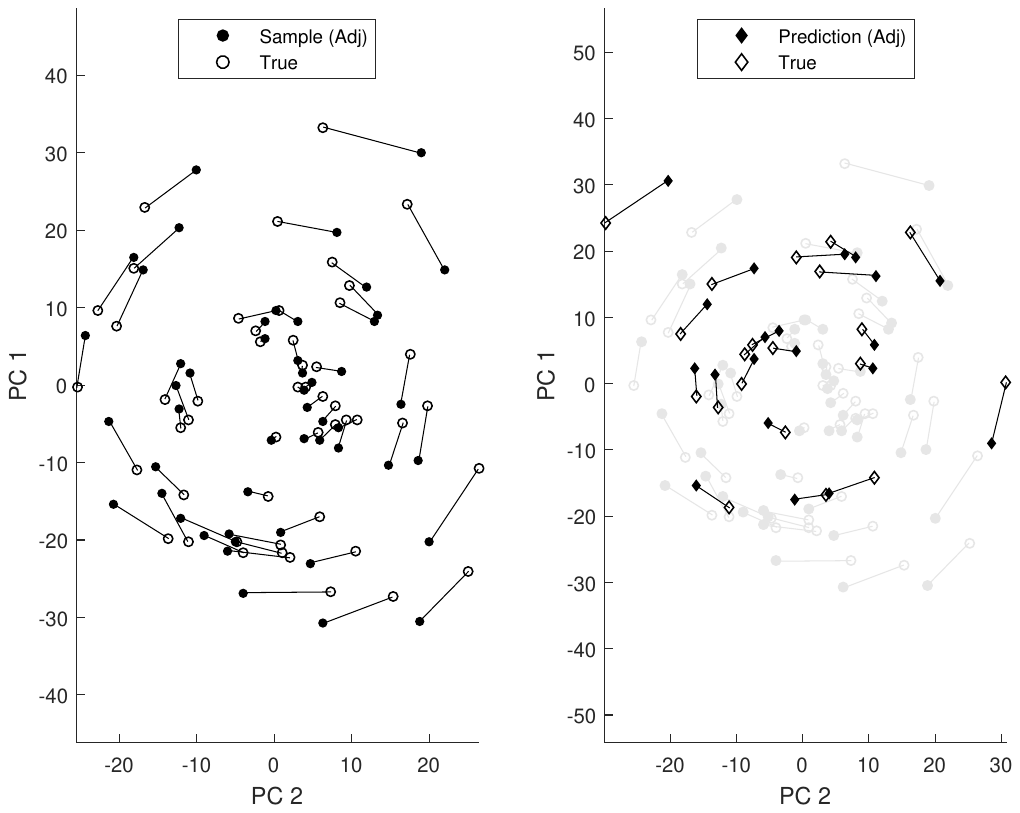}
  \vskip -4cm
  \caption{Bias-adjusted sample and prediction scores using (\ref{eq:bias-adjustment}) for the toy data introduced in Fig.~\ref{fig:projectionScores2}. The estimates (\ref{eq:rho-est1}) are $(\tilde\rho_1,\tilde\rho_2) = (1.385,1.546)$, very close to the theoretical values $(\rho_1,\rho_2) = (1.385,1.557)$.  Both sample and prediction scores are simultaneously rotated about 16 degrees clockwise.
  \label{fig:projectionScores2-repeat}}
\end{figure}

An application of the above bias-adjustment procedure is exemplified in Fig.~\ref{fig:projectionScores2-repeat}. There, the magnitudes of the sample and prediction scores are well-adjusted. 

Our next proposed estimators are motivated by the well-known jackknife bias adjustment procedures and also by the leave-one-out cross-validation. For simplicity, assume $m=1$. The bias-adjustment factor we aim to estimate is $\rho_1 = ({1 + \tau^2/ \| \xi_1 \|_2^2})^{1/2}$, where $\xi_1 = d^{-1/2}w_1 = \sigma_1(z_{11},\ldots, z_{1n})\T$ is the scaled true scores for the first principal component.

Write, for each $j = 1,\ldots,n$, the $j$th scaled sample score as $\hat{\varpi}_{1j} = d^{-1/2}\hat{u}_1\T X_j$, and the $j$th scaled prediction score as
$$\hat{\varpi}_{1(j)} = d^{-1/2}\hat{u}_{1(-j)}\T X_j,$$
where $\hat{u}_{1(-j)}$ is the first principal component direction, computed from $\Xc_{(-j)}$, i.e., the data except the $j$th observation.

From Theorem~\ref{thm:main}, $\rho_1$ is the asymptotic bias-adjustment factor for $\hat{\varpi}_{1}$; $\hat{\varpi}_{1j} = \rho_1 \varpi_{1j} + O_p(d^{-1/4})$.
For $\hat{\varpi}_{1(j)}$, again applying Theorem~\ref{thm:main}, we get
$\hat{\varpi}_{1(j)} =  \rho_{1(-j)}^{-1} {\varpi}_{1j} + O_p(d^{-1/2}),$
 where $\rho_{1(-j)} = ({1 + \tau^2/ \| \varpi_{1(-j)} \|_2^2})^{1/2} $ is the bias-adjustment factor computed from $\Xc_{(-j)}$, using
$\varpi_{1(-j)} =  \sigma_1(z_{11},\ldots,z_{1,{j-1}},z_{1,j+1},\ldots z_{1n})\T$.
To simplify terms, Taylor expansion is used to expand
$\rho_{1(-j)}$ as a function of $\varpi^2_{1j}/n$, resulting in
\begin{equation}\label{eq:jack3}
\rho_{1(-j)} = \left(1 + \frac{\tau^2 / n}{ \| \varpi_{1} \|_2^2/n - \varpi_{1j}^2/n } \right)^{1/2} = \rho_1 + \frac{1}{2\rho_1 } \frac{\| \varpi_{1} \|_2^2/n}{\tau^2}{} \frac{\varpi^2_{1j}}{n} + O_p(\frac{1}{n^2}).
\end{equation}

Using the approximation $$\rho_1\rho_{1(-j)} \approx \rho^2_1 + \frac{\| \varpi_{1} \|_2^2}{2\tau^2}{} \frac{\varpi^2_{1j}}{n^2}$$ given by (\ref{eq:jack3}), we write the  ratio of the sample and prediction scores to cancel out the unknown true score $\varpi_{1j}$ as follows:
$$\left(\frac{\hat{w}_{1j}}{\hat{w}_{1(j)}}\right)^{1/2}  =  \left(\frac{\hat{\varpi}_{1j}}{\hat{\varpi}_{1(j)}}\right)^{1/2}
      \approx \rho_1.$$
Based on the above heuristic, we define the following estimators of the bias-adjustment factors:
\begin{align}
\hat{\rho}_i^{(1)}
&= \frac{1}{n}\sum_{j=1}^n \left(\frac{\hat{w}_{ij}}{\hat{w}_{i(j)}}\right)^{1/2} , \label{eq:rho-est2} \\
\hat{\rho}_i^{(2)}  &=  \left( \frac{ \sum_{j=1}^n \hat{w}_{ij}  } { \sum_{j=1}^n \hat{w}_{i(j)}} \right)^{1/2} , \label{eq:rho-est3}\\
\hat{\rho}_i^{(3)}  &=  \left( \frac{ \sum_{j=1}^n \hat{w}^2_{ij}  } { \sum_{j=1}^n \hat{w}^2_{i(j)}} \right)^{1/4}. \label{eq:rho-est4}
\end{align}
In implementing the above estimators, we used absolute values of the sample and predicted scores. The estimator (\ref{eq:rho-est4}) is a ratio of the sample and prediction score variances, obtained by a leave-one-out estimation of prediction scores.

The estimators $\hat{\rho}_i^{(1)}$, $\hat{\rho}_i^{(2)}$, and $\hat{\rho}_i^{(3)}$ tend to overestimate $\rho$ for small sample size $n$, as expected from (\ref{eq:jack3}). In our numerical experiments, these three estimators perform similarly.

\section{Numerical studies}\label{sec:numericalstudies}

\subsection{Simulations to confirm the asymptotic bias and near-perfect correlations}\label{sec:numeric-theory}
In this section, we compare the theoretical asymptotic quantities derived in Section~\ref{sec:mainresults} with their finite-dimensional empirical counterparts.

\begin{figure}[tp]
  \centering
  \vskip -3cm
  \includegraphics[width=0.9\textwidth]{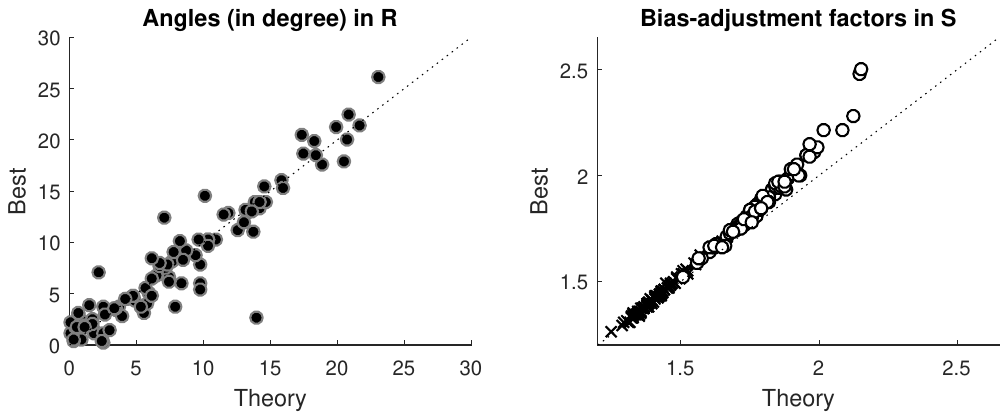}
  \vskip -8cm
  \vskip -3cm
  \includegraphics[width=0.9\textwidth]{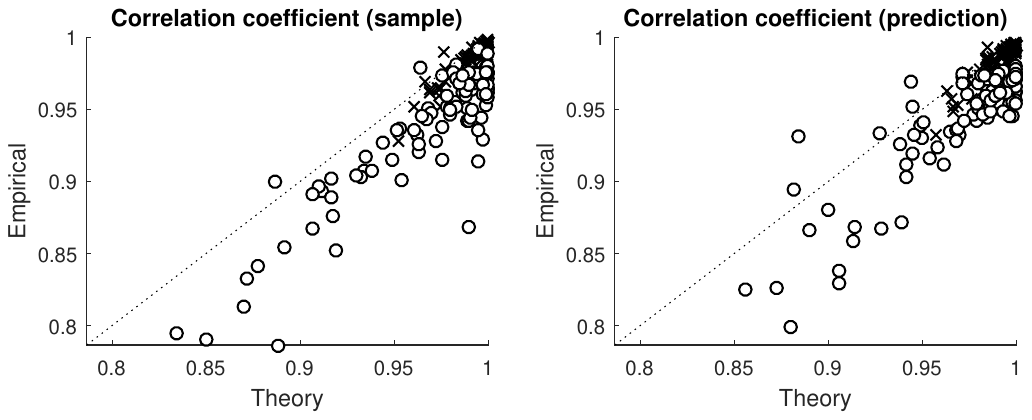}
  \vskip -8cm
  \caption{(Top row) Theoretical rotation angles $\theta_R$ and bias-adjustment factors $\rho_1$ ($\times$), $\rho_2$ ($\bigcirc$), compared with the best-fitting Procrustes counterparts ($\check\theta_R, \check{\rho}_i(W_1)$). (Bottom row) Empirical correlations compared with their limits in Theorem~\ref{thm:correlations}.
  \label{fig:theory1}}
\end{figure}

First, the theoretical values of the scaling bias $\rho_i$ and the  rotation matrix $R$ in Theorem~\ref{thm:main} are compared with their empirical counterparts. The empirical counterparts of the two matrices $R, S$ are defined as the minimizer of the Procrustes problem
\begin{equation}
\label{eq:procrustes}
\min \norm{W_1 - \widehat{W}_1\T S_0^{-1} R_0}_F^2,
\end{equation}
with the constraint that $S_0$ is a diagonal matrix with positive entries and $R_0$ is an orthogonal matrix. The solutions are denoted by $\widecheck{S} = \mbox{diag}(\check{\rho}_1(W_1),\ldots, \check{\rho}_m(W_1))$ and $\widecheck{R}$. For simplicity, we consider the $m = 2$ case, and parameterize $R$ by the rotation angle, $\theta_R = \cos^{-1}(R_{1,1})$, and $\widecheck{R}$ by $\check\theta_R = \cos^{-1}(\widecheck{R}_{1,1})$.
We compare $\theta_R$ with $\check\theta_R$ and $\rho_i(W_1)$ with $\check\rho_i(W_1)$, from a 2-component model with $(n,d) = (50, 5000)$ (precisely, the spike model with $m=2$ and $\beta = 0.3$ in Section \ref{sec:SIMmodels}).
Note that both the theoretical values and the best-fitted values depend on the true scores $W_1$. To capture the natural variation given by $W_1$, the experiment is repeated for   100 times.
The results, summarized in the top row of Fig.~\ref{fig:theory1}, confirm that the asymptotic statements in Theorem~\ref{thm:main} approximately hold for finite dimensions. In particular, the rotation matrices $R$ and $\widecheck{R}$ are very close to each other. The Procrustes-fitted, or ``best'', $\check\rho_i$ tends to be larger than the asymptotic, or theoretical, $\rho_i$, especially for $i = 2$ (shown as $\bigcirc$ in Fig.~\ref{fig:theory1}) and for larger values of $\rho_2$. This is not unexpected. Larger values of $\rho_2$ are from smaller $\lambda_2(\Wc)$. Take an extreme case where $\lambda_2(\Wc) = 0$, then by (\ref{eq:score-sample-diverge}) in  Theorem~\ref{thm:main}, the sample scores are of magnitude $d^{1/2}$ compared to the true scores. Thus, as $\lambda_2(\Wc)$ decreases to 0, the Procrustes scaler $\check\rho_2$ empirically interpolates the finite-scaling case (\ref{eq:scarot1}) to the diverging case (\ref{eq:score-sample-diverge})  of  Theorem~\ref{thm:main}.

Second, we compare the limit of correlation coefficients in Theorem~\ref{thm:correlations} with finite-dimensional empirical correlations, $r(\hat{w}_k,w_k)$, for $k = 1, 2$. For the correlation coefficient of the prediction scores, we use the sample correlation coefficient between $(\hat{w}_{k*},w_{k*})$, as an estimate of $\mbox{Corr}(\hat{w}_{k*},w_{k*} \mid W_1)$.
The simulated results are shown in the bottom row of Fig.~\ref{fig:theory1}.
The empirical correlation coefficients tend to be smaller than the theoretical counterparts, but both are higher for stronger ``signal strength" $n\sigma_k^2 = \E(\lambda_k(\Wc))$.

Third, from the same simulations, it can be checked that the $k$th, where $k>m$, sample scores are diverging, while the prediction scores are stable, as indicated in (\ref{eq:score-sample-diverge}) and (\ref{eq:score-prediction-stable}). To confirm this, we choose $k = 3$ and for each experiment, compute $\widehat{\Var}(\hat{w}_3)$, the sample variance of the sample scores, and an approximation of  ${\Var}(\hat{w}_{3*})$. The results are shown in Table~\ref{tab:k=3case}. As expected, the sample scores are grossly inflated, while the prediction scores are stable. Finally, the conjecture in Remark~\ref{remark:conjecture} is also empirically checked; Table~\ref{tab:k=3case} also shows that for large $d$, the sample (or prediction) and true scores for the $k$th, $k>m$, component are nearly uncorrelated.

\begin{table}[t]
\centering
\begin{tabular}{c|cc}
          & Sample scores & Prediction scores  \\
  \hline
 Variance & 120.7(4.4) & 1.38(0.2) \\
 Corr. Coef. & -0.0024(0.2) & -0.004(0.15)
\end{tabular}
\caption{The $k$th sample and prediction scores (unadjusted) for the case $k>m$. Shown are the mean (standard deviation) of the variances and correlation coefficients to true scores, from 100 repetitions. The true variance is $\lambda_3 = {\Var}(w_{3*}) \approx 6.5$.
  \label{tab:k=3case}}
\end{table}

\subsection{Numerical performance of the bias-adjustment factor estimation}\label{sec:SIMmodels}

We now test our estimators of the bias-adjustment factor $\rho_i$, using the following data-generating models with $m = 2$.

The first one is called a \emph{spike model}. We sample from the $d$-dimensional zero-mean normal distribution where the first two largest eigenvalues of the covariance matrix are
$\lambda_i = \sigma_i^2 d$, for $i = 1,2$, where $(\sigma_1^2,\sigma_2^2) = (0.02, 0.01)$. The rest of eigenvalues are slowly-decreasing. In particular, $\lambda_i = \tau i^{-\beta}$, where $\tau = [ \sum_{i=3}^d i^{-\beta}/(d-2) ]^{-1}$. We set $\beta = 0.3$ or $0.5$. This spike model has more than two unique principal components for each fixed dimension, but in the limit $d\to\infty$, only the first two principal components are useful.

The second model is a \emph{mixture model}. Let $\mu_g$ ($g=1,2,3$) be $d$-dimensional vectors, the elements of which are randomly drawn from $\{-a,0,a\}$ with replacement for a given $a>0$, then assumed as fixed quantities. Given $\mu_g$'s we sample from the mixture model $X \mid G = g \sim N(\mu_g, \Id_d)$, $P(G=g) = p_g>0$, $\sum_{g=1}^3 p_g = 1$. We set $(p_1,p_2,p_3) = (0.5,0.3,0.2)$. 
It can be checked that $\Cov(X)$  satisfies the assumption of the 2-component model in (A1)--(A4). 

For various cases of high-dimension, low-sample-size situations, ranging $d =5,000$ to $20,000$ and $n = 50$ to $100$, random samples from each of these models are generated.
For each case, the theoretical quantity $\rho_i = \rho_i(W_1)$ and the best-fitted Procrustes scaler $\check\rho_i = \check\rho_i(W_1)$ are computed. These quantities depend on the $m\times n$ random matrix $W_1$. The mean and standard deviation of $\rho_i$ (from 100 repetitions) are shown in the first column of  Table~\ref{tab:simulation}. As expected, the theoretical value $\rho_i$ depends on the sample size $n$; large sample size decreases the bias, $\E(\rho_i)$, and also decreases the variance $\Var(\rho_i)$. 

The mean of the best-fitted scaler $\check\rho_i$ ($i=1$) is displayed in the second column of the table. While they are quite close to the theoretical counterpart, $\check\rho_i$s are significantly larger for the mixture model, whose signal-to-noise ratio is smaller than the spike model, and for the not-so-large dimension $d = 5,000$. This is not unexpected, since the theoretical values are also based on the dimension-increasing asymptotic arguments.

We further compute the proposed estimators of $\rho_i$, given in (\ref{eq:rho-est1}), (\ref{eq:rho-est2})--(\ref{eq:rho-est4}). We also compute the estimator derived from  \cite{lee2010convergence}, which is the square-root of the reciprocal of the shrinkage factor, obtained by numerical iterations, denoted by $\hat{d}_\nu$ in \cite{lee2010convergence}. (The relation of \cite{lee2010convergence} to our work is further discussed in Section~\ref{sec:discussion}.) 
 All of the methods considered provide accurate estimates of the theoretical quantity $\rho_i$. We omit the numerical results from the estimators (\ref{eq:rho-est3}) and (\ref{eq:rho-est4}), as their performances are very close to those from (\ref{eq:rho-est2}). The supplementary material contains an extended table of Table~\ref{tab:simulation}, including the case for $\rho_2$.

\begin{table}[tp]
{\footnotesize
\centering
\begin{tabular}{cccccccc}
& & & & & $\rho_1$ \\
\cline{4-8}
   & $d$     & $n$   & Theory & Best & Asymp. & Jackknife & LZW \\
       \hline
 & 5000  & 50  & 1.41 (0.07) & 1.42 & 1.40 & 1.43 & 1.41 \\
Spike model & 10000 & 50  & 1.42 (0.06) & 1.43 & 1.42 & 1.44 & 1.42 \\
$\beta = 0.3$ & 10000 & 100 & 1.23 (0.03) & 1.23 & 1.23 & 1.24 & 1.23 \\
 & 20000 & 100 & 1.23 (0.02) & 1.23 & 1.23 & 1.24 & 1.23 \\
 \hline
 & 5000  & 50  & 1.42 (0.08) & 1.45 & 1.41 & 1.45 & 1.40 \\
Spike model & 10000 & 50  & 1.43 (0.07) & 1.45 & 1.43 & 1.46 & 1.42 \\
$\beta = 0.5$ & 10000 & 100 & 1.22 (0.02) & 1.23 & 1.22 & 1.23 & 1.21 \\
 & 20000 & 100 & 1.23 (0.02) & 1.23 & 1.23 & 1.24 & 1.22 \\
 \hline
& 5000  & 50  & 2.06 (0.06) & 2.22 & 1.92 & 2.14 & 2.00 \\
Mixture model & 10000 & 50  & 2.09 (0.06) & 2.17 & 1.98 & 2.14 & 2.02 \\
$a = 0.15$ & 10000 & 100 & 1.63 (0.02) & 1.67 & 1.61 & 1.65 & 1.63 \\
 & 20000 & 100 & 1.64 (0.02) & 1.66 & 1.62 & 1.66 & 1.63
\end{tabular}%
}
\caption{Simulation results from 100 repetitions. ``Theory'' is mean (standard deviation) of $\rho_i$; ``Best'' is $\check\rho_i$ (\ref{eq:procrustes}); ``Asymp.'' is $\tilde\rho_i$ (\ref{eq:rho-est1}); ``Jackknife'' is $\hat\rho_i^{(1)}$ (\ref{eq:rho-est2}); ``LZW'' is from \cite{lee2010convergence}. Averages are shown for the latter four columns. The standard errors of the quantities in estimation of $\rho_i$ are at most 0.04.
  \label{tab:simulation}}
\end{table}

\subsection{Bias-adjustment improves classification}
Our last simulation study is an application of the bias-adjustment procedure to classification. Our training and testing data, each with sample size $100$, are sampled from the mixture model with three groups, as described in Section~\ref{sec:SIMmodels}.
As frequently used in practice \citep{adam2008classification}, dimension reduction by the standard principal component analysis is performed first, then a classification rule by the support vector machine \citep[SVM, ][]{Cristianini2000} is trained on the sample principal component scores. In this simulation, we fix $m =2$ and $d = 5000$.
We compare the training and testing missclassification error rates (estimated by 100 repetitions) of the SVMs trained (and tested) either on the unadjusted sample and prediction scores, $\widehat{W}_1$ and $\widehat{W}_\star$, or on the bias-adjusted sample and prediction scores, $\widehat{W}_1^{\rm (adj)}$ and $\widehat{W}_\star^{\rm (adj)}$ in (\ref{eq:bias-adjustment}). The estimated error rates are shown in Table~\ref{tab:classification}. It is clear that the use of bias-adjusted scores greatly improves the performance of  classification.

\begin{figure}[t]
  \centering
  \vskip -3cm
  \includegraphics[width=1\textwidth]{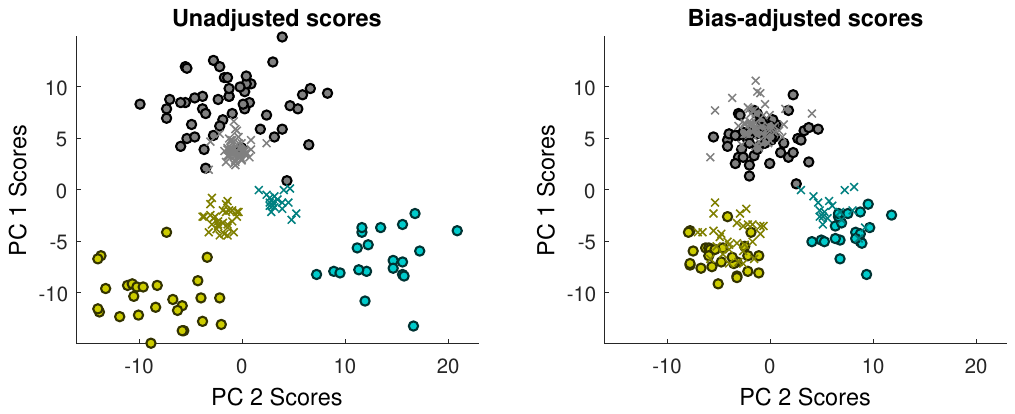}
  \vskip -8cm
  \caption{Bias-adjusted scores from the mixture models greatly improve the classification performance. Different colors correspond to different groups. Symbol $\bigcirc$ represents the sample scores (unadjusted in the left, adjusted in the right); symbol $\times$ represents the prediction scores. \label{fig:sim-classify}}
\end{figure}

To better understand the huge improvement of classification performances, we plot the sample and prediction scores that are inputs of the classifier. In  Fig.~\ref{fig:sim-classify}, the classifier is estimated from the the sample scores (symbol $\bigcirc$) and is used to classify future observations, i.e. the prediction scores (symbol $\times$). Due to the scaling bias, the unadjusted sample and prediction scores are of different scales (shown in the left panel), and classification is bound to fail. On the other hand, the proposed bias-adjustment, shown in the right panel, works well for this data set, leading to a better classification performance.

\begin{table}[t]
\centering
\begin{tabular}{c|cc}
          & Unadjusted scores & Bias-adjusted scores  \\
  \hline
 Training Error & 0.04(0.02) & 0.07(0.03) \\
 Testing Error &  21.4(1.33) & 1.98(0.23)
\end{tabular}
\caption{Means (standard errors) of Missclassification error rates (in percent).
  \label{tab:classification}}
\end{table}

\section{Discussion}\label{sec:discussion}

The standard principal component analysis is shown to be useful in the dimension reduction of data from the $m$-component models with diverging variances. In particular, in the high-dimension, low-sample-size asymptotic scenario we reveal that the sample and prediction scores have systematic biases that can be consistently adjusted.  We propose several estimators of the scaling bias, while there is no compelling reason to adjust rotational bias. The amount of bias is large when the sample size is small and when the variance of accumulated noise is large compared to the variances of the first $m$ components.

\cite{lee2010convergence} discussed adjusting bias in the prediction of principal components, based on the random matrix theory and the asymptotic scenario of $d/n \to \gamma \in (0,\infty)$, $n \to \infty$. They showed that the prediction scores tend to be smaller than the sample scores, and the ratio of the shrinkage is asymptotically ${\rm sd}(\hat{w}_{i1}) / {\rm sd}(\hat{w}_{i*}) \approx \rho_i^{\rm (LZW)} = \frac{\lambda_i - 1}{\lambda_i + \gamma -1}$. This ``shrinkage factor"  $\rho_i^{\rm (LZW)}$ corresponds to the squared reciprocal of our scaling bias, $\rho_i^{-2}$. 
Our work can be thought of as an extension of \cite{lee2010convergence} from the asymptotic regime $d \asymp n$ to the high-dimension, low-sample-size situations (see also \cite{lee2014convergence,dey2019asymptotic}).
Finally, we note that in the asymptotic scenario of \cite{lee2010convergence,lee2014convergence} and \cite{dey2019asymptotic} there is no rotational bias. This is because in their limit the sample size is infinite. We show that the rotational bias is universal to both sample and prediction scores and is of order $n^{-1/2}$.

\section*{Supplementary Materials}

This supplementary material contains technical details and proofs, and a table summarizing simulation results.

\appendix

\setcounter{section}{0}
\setcounter{equation}{0}
\def\theequation{S\arabic{section}.\arabic{equation}}
\def\thesection{S\arabic{section}}


\section{Proofs of Theorem 1 and Lemma 2}

For reference, we restate the theorems and formulas in the main article that are used in the proof.

\begin{theorem} \label{thm:main}
Assume the $m$-component model under Conditions (A1)--(A4) and let $n > m \ge 0$ be fixed and $d\to\infty$. Then, the first $m$ sample and prediction scores are systematically biased:
    \begin{align}
     \widehat{W}_1 &= S R\T W_1  + O_p(d^{-1/4}), \label{eq:scarot1}\\
     \widehat{W}_{*} &= S^{-1} R\T  W_*  + O_p(d^{-1/2}),  \label{eq:scarot2}
     \end{align}
    where $R = [v_1(\Wc),\ldots, v_m(\Wc)]$, $S = \mbox{diag}(\rho_1,\ldots,\rho_m),$  
and $\rho_k = \sqrt{ 1+ \tau^2/\varphi_k(\Wc) }$.
 Moreover,
 for $k > m$,
     \begin{align}
     \hat{w}_{kj} &= O_p(d^{1/2}), \quad j = 1,\ldots,n, \label{eq:score-sample-diverge}\\
      \hat{w}_{k*} &= O_p(1).  \label{eq:score-prediction-stable}
     \end{align}
\end{theorem}

\begin{lemma}\label{thm:PCA-asymptotic} [Theorem S2.1, \cite{Jung2017}]
 Assume the conditions of Theorem~\ref{thm:main}.
  (i)  the sample principal component variances converge in probability as $d \to \infty$;
   \begin{equation*}
  d^{-1} n \hat\lambda_{i} = \left\{
                                   \begin{array}{ll}
                                     \varphi_i (\Wc) + \tau^2 + O_p(d^{-1/2}) , & i = 1,\ldots,m; \\
                                     \tau^2 + O_p(d^{-1/2}), & i = m+1,\ldots,n.
                                   \end{array}
                                 \right.
   \end{equation*}
  (ii) The inner product between sample and population PC directions converges  in probability as $d \to \infty$;
   \begin{equation*}
   \hat{u}_{i}\T u_{j} = \left\{
                                     \begin{array}{ll}
                                       \rho_i^{-1} v_{ij}(\Wc) + O_p(d^{-1/2}) , & i,j = 1,\ldots,m; \\
                                       O_p(d^{-1/2}), & \mbox{ otherwise.}
                                     \end{array}
                                   \right.
   \end{equation*}
\end{lemma}

\begin{lemma}\label{lem:error_term_bound}
 Assume the $m$-component model with (A1)--(A4) and let $n > m \ge 0$ be fixed. For $k = 1,\ldots,n$, $\E( \epsilon_{k*} | W_1) = 0$, and
   \begin{align} \label{eq:variancelimit}
    \lim_{d\to\infty}\Var ( \epsilon_{k*}   \mid W_1)
    &= \upsilon^2_O  / (\varphi_k(\Wc) + \tau^2 ), \quad \hbox{for $k \le m$;} \\
\lim_{d\to\infty} \frac{1}{n-m} \sum_{k = m+1}^n \Var ( \epsilon_{k*}   \mid W_1)
    &= \upsilon^2_O  / \tau^2   ,  \label{eq:variancelimit2}
\end{align}
where $\upsilon^2_O = \lim_{d\to\infty} d^{-1}\sum_{i=m+1}^d \lambda^2_i$. As $d\to\infty$,  $\epsilon_{k*}  = O_p(1)$.
\end{lemma}

\begin{proof}[Proof of Lemma~\ref{lem:error_term_bound}]
Fix $k = 1,\ldots,n$.
Let $Y_i = \sqrt{\lambda_i}z_{i*} p_{ki}$, where $p_{ki} = \hat{u}_k\T u_i$. Then $\epsilon_{k*} =  \sum_{i=m+1}^d Y_i$.
Since $z_{i*}$ and $p_{ki}$ are independent, for each $i > m$, $\E(Y_i \mid W_1) = 0$ and
\begin{align*}
\Var (\sum_{i=m+1}^d Y_i  \mid W_1)  = \E(\sum_{i=m+1}^d {\lambda_i}z^2_{i*} p_{ki}^2  \mid W_1)  =
 \sum_{i=m+1}^d {\lambda_i}  \E(p_{ki}^2  \mid W_1) ,
\end{align*}
where we use the fact that $\E(z_{i*}) = 0$, $\E(z_{i*}^2) = 1$.

For $k \le m$, if the following claim,
\begin{equation} \label{eq:proof:Epki2}
\E(p_{ki}^2 \mid W_1 )
 = d^{-1} \frac{\lambda_i}{(\lambda_k(\Wc) + \tau^2)} + O(d^{-3/2}),
\end{equation}
is true for any $i > m$, then it is easy to check (\ref{eq:variancelimit}).

To show (\ref{eq:proof:Epki2}), we first post-multiply $\hat{v}_i$ to
\begin{equation} \label{eq:SVD on X}
\Xc = \sqrt{n} \sum_{i=1}^n  \sqrt{\hat\lambda_i} \hat{u}_i \hat{v}_i\T,
\end{equation}
to obtain $\hat{u}_i = (n\hat\lambda_i)^{-1/2} \Xc \hat{v}_i$. By writing $z_i\T = \lambda_i^{-1/2} w_i\T = (z_{i1},\ldots,z_{in})$, we have
\begin{align*}
p_{ki} & = u_i\T \hat{u}_k \\
       & = (n\hat\lambda_k)^{-1/2} u_i\T \Xc \hat{v}_k \\
       & = (n\hat\lambda_k)^{-1/2} \lambda_i^{1/2} z_i\T \hat{v}_k.
\end{align*}
Thus,
\begin{align}
p_{ki}^2 & = d^{-1}\frac{\lambda_i}{nd^{-1}\hat\lambda_k}    ( z_i\T \hat{v}_k )^2 \nonumber \\
         & = d^{-1}\frac{\lambda_i}{\lambda_k(\Wc) + \tau^2 + O_p(d^{-1/2})}   ( z_i\T \hat{v}_k )^2  \nonumber \\
         & = d^{-1} \frac{\lambda_i}{\lambda_k(\Wc) + \tau^2} ( z_i\T \hat{v}_k )^2 + O_p(d^{-3/2}). \label{eq:proof-pki2}
\end{align}
In (\ref{eq:proof-pki2}), we used Lemma~\ref{thm:PCA-asymptotic}(i) and that $(1+x)^{-1} = 1 + O(x)$, and the fact that $|z_i\T \hat{v}_k|^2 \le \norm{z_i}^2_2 \norm{\hat{v}_k}_2^2 = \norm{z_i}_2^2 = O_p(1)$.

Write $( z_i\T \hat{v}_k )^2 =  [ z_i\T {v}_k(W_1\T W_1) + z_i\T (\hat{v}_k - v_k(W_1\T W_1)) ]^2$. Note that $W_1\T W_1$ is an $n \times n$ matrix, and is different from the $m \times m$ matrix $\Wc = W_1 W_1\T$.
It can be shown that the right singular vector $\hat{v}_k$ converges in probability to $v_k(W_1\T W_1)$ \citep[see, e.g., Lemma S1.1 of][]{Jung2017}: For $k =1,\ldots,m$,
\begin{equation}
\label{eq:v_to_v(W)}
\hat{v}_k = v_k(W_1\T W_1) + O_p(d^{-1/2}).
\end{equation}
Thus we get $|z_i\T (\hat{v}_k - v_k(W_1\T W_1))| \le \norm{z_i}_2 \norm{\hat{v}_k - v_k(W_1\T W_1))}_2 = O_p(d^{-1/2})$. Therefore,
\begin{align}
\E( ( z_i\T \hat{v}_k )^2  \mid W_1)
&=   \E( ( z_i\T {v}_k(W_1\T W_1))^2  \mid W_1) + O(d^{-1/2}) \nonumber \\
& =  \sum_{\ell = 1}^n  \E( z_{i\ell}^2) v^2_{k\ell}(W_1\T W_1) + O(d^{-1/2})\nonumber  \\
&= 1 + O(d^{-1/2}). \label{eq:proof-pki2_2}
\end{align}
Combing (\ref{eq:proof-pki2}) and (\ref{eq:proof-pki2_2}), we get (\ref{eq:proof:Epki2}) for $k\le m$ as desired.

To show (\ref{eq:variancelimit2}), note that $\Wc = W_1 W_1\T$ is of rank $m$.
For $k> m$, with $\lambda_k(\Wc)=0$, (\ref{eq:proof-pki2}) holds. Thus,
\begin{align}
\frac{1}{n-m}\sum_{k=m+1}^n \Var(\epsilon_{k*} \mid W_1)
  & =  \frac{1}{n-m}\sum_{k=m+1}^n \sum_{i=m+1}^d \lambda_i \E(p_{ki}^2 \mid W_1) \label{eq:proof:last2} \\
  & = \frac{1}{d(n-m)}\sum_{i=m+1}^d  \lambda_i^2 / \tau^2 \sum_{k=m+1}^n \E( (z_i\T \hat{v}_k)^2 \mid W_1). \nonumber
\end{align}
To simplify the expression $\E( (z_i\T \hat{v}_k)^2 \mid W_1)$, one should not naively try (\ref{eq:proof-pki2_2}). This is because that (\ref{eq:proof-pki2_2}) does not apply for $k>m$ due to the non-unique $k$th eigenvector $v_k(W_1\T W_1)$ of the rank-$m$ matrix $W_1\T W_1$. Instead, from
$$
\sum_{k=m+1}^n (z_i\T \hat{v}_k)^2 = z_i\T z_i - \sum_{k=1}^m (z_i\T \hat{v}_k)^2,
$$
and (\ref{eq:proof-pki2_2}) for $k\le m$, we get
\begin{equation}\label{eq:proof:last}
\sum_{k=m+1}^n \E( (z_i\T \hat{v}_k)^2 \mid W_1) = n - m + O(d^{-1/2}).
\end{equation}
Taking the limit $d\to\infty$ to (\ref{eq:proof:last2}), combined with (\ref{eq:proof:last}), leads to (\ref{eq:variancelimit2}).

The last statement, $\epsilon_{k*} = O_p(1)$, easily follows from the fact
$\lim_{d\to\infty}\Var(\epsilon_{k*}) \le \upsilon_O^2 / \tau^2 (n-m) <\infty$, which is obtained by (\ref{eq:variancelimit}) and (\ref{eq:variancelimit2}).
\end{proof}

We are now ready to show Theorem~\ref{thm:main}.
Note that the results on the sample scores, (\ref{eq:scarot1}) and (\ref{eq:score-sample-diverge}), can be easily shown, using the decomposition
$d^{-1/2}\hat{w}_k = \sqrt{d^{-1}n\hat\lambda_k}\hat{v}_k$, together with Lemma~\ref{thm:PCA-asymptotic}(i) and (\ref{eq:v_to_v(W)}).
We show (\ref{eq:scarot2}) and (\ref{eq:score-prediction-stable}).

\begin{proof}[Proof of Theorem~\ref{thm:main}]
Proof of (\ref{eq:scarot2}).
Recall the decomposition
\begin{align}\label{eq:projectionNew}
\hat{w}_{k*} = \hat{u}_k\T X_{*} = \sum_{i=1}^m w_{i*} \hat{u}_k\T u_i + \epsilon_{k*},
\end{align}
where $\epsilon_{k*} =  \sum_{i=m+1}^d w_{i*} \hat{u}_k\T  u_i.$
 Using the notation $p_{ki} = \hat{u}_k\T  u_i$, we write
$\hat{w}_{k*} = (p_{k1},\ldots,p_{km} ) (w_{1*},\ldots,w_{m*})\T + \epsilon_{k*}$. Putting all parts together, we have
$$
\widehat{W}_*= d^{-1/2}(\hat{w}_{1*},\ldots,\hat{w}_{m*})\T =   \left(  \begin{array}{ccc}
                                           p_{11} & \cdots & p_{1m} \\
                                           \vdots & \ddots & \vdots \\
                                           p_{m1} & \cdots & p_{mm} \\
                                         \end{array}
                                       \right) {W}_*   + \tilde\epsilonv_{k*},
$$
where $\tilde\epsilonv_{k*} = d^{-1/2}( \epsilon_{1*},\ldots, \epsilon_{m*})\T$. By Lemma~\ref{lem:error_term_bound}, as $d\to\infty$, $\tilde\epsilonv_{k*} = O_p(d^{-1/2})$.
 Since $p_{ki} = \rho_{k}^{-1} v_{ki}(\Wc) + O_p(d^{-1/2})$, by Lemma~\ref{thm:PCA-asymptotic}(ii), we have
$$ \widehat{W}_*\T=  S^{-1} R\T {W}_*\T    + O_p(d^{-1/2}).$$

Proof of (\ref{eq:score-prediction-stable}). Using the decomposition (\ref{eq:projectionNew}), and by the fact
$\epsilon_{k*} = O_p(1)$,
from Lemma~\ref{lem:error_term_bound},
it is enough to show $\sum_{i=1}^m w_{i*} p_{ki} = O_p(1)$.
But, since Lemma~\ref{thm:PCA-asymptotic} implies $d^{\half} p_{ki} = O_p(1)$ for any pair of $(k,i)$ such that $k>m, i\le m$, we have
$
\sum_{i=1}^m w_{i*} p_{ki} = \sigma_i (d^{\half}p_{k1},\ldots,d^{\half} p_{km})(z_{1*},\ldots,z_{m*}) = O_p(1),
$
\end{proof}

\section{Proof of Theorem 2}
\setcounter{equation}{0}

\begin{theorem}\label{thm:correlations}
 Let $\zeta_{kj} = \lambda_k(\Wc)/(\sum_{\ell= 1}^m v^2_{\ell j}(\Wc) \lambda_\ell(\Wc) )$ and $\bar\zeta_{kj}  = \sigma^2_k/(\sum_{\ell= 1}^m v^2_{\ell j}(\Wc) \sigma^2_\ell )$.
Under the assumptions of Theorem~\ref{thm:main}, as $d\to\infty$,  for  $k,j = 1\ldots,m$,
\begin{itemize}
\item[(i)] $r(\hat{w}_k,w_j) \to v_{kj}(\Wc) \zeta_{kj}^{1/2}$ in probability ;
\item[(ii)] $ \lim_{d\to\infty}\mbox{Corr}(\hat{w}_{k*},w_{j*} \mid W_1) = v_{kj}(\Wc) \bar\zeta_{kj}{ }^{1/2}$.
\end{itemize}
\end{theorem}

\begin{proof}[Proof of Theorem~\ref{thm:correlations}]
Proof of (i). Write the singular value decomposition of the $m \times n $ matrix of scaled scores $W_1$ as
\begin{equation}\label{eq:W_1-svd}
W_1 = R \mbox{diag}(\sqrt{\lambda_1(\Wc)},\ldots, \sqrt{\lambda_1(\Wc)}) G\T,
\end{equation}
where $G = [g_1,\ldots, g_m]$ is the $n \times m $ matrix consisting of right singular vectors of $W_1$. The left singular vector matrix $R = [v_1(\Wc),\ldots,v_m(\Wc)]$ is exactly the matrix $R$ appearing in Theorem~\ref{thm:main}. Since
$$ W_1 = \sum_{\ell = 1}^m \sqrt{\lambda_\ell(\Wc)} v_{\ell}(\Wc)g_{\ell}\T,$$
the $j$th row of $W_1$ is, for $j \le m$,
$$d^{-\half}w_j\T = \sum_{\ell = 1}^m \sqrt{\lambda_\ell(\Wc)} v_{\ell j}(\Wc)g_{\ell}\T.$$

For the scaled sample score $d^{-1/2}\hat{w}_k$, $k\le m$, we obtain from Theorem~\ref{thm:main} and (\ref{eq:W_1-svd}) that
$\widehat{W}_1 = S \mbox{diag}(\sqrt{\lambda_1(\Wc)},\ldots, \sqrt{\lambda_1(\Wc)}) G\T + O_p(d^{-1/4})$
and its $k$th row $d^{-1/2}\hat{w}_k = \sqrt{\lambda_k(\Wc) + \tau^2} g_k   + O_p(d^{-1/4})$.
Since $g_\ell$'s are orthonormal,
\begin{equation*}
\| d^{-\half}\hat{w}_k\|_2 = \sqrt{\lambda_k(\Wc) + \tau^2} + O_p(d^{-1/4}),
\end{equation*}
and
\begin{align*}
d^{-1} \hat{w}_k\T w_j &= (d^{-1/2}\hat{w}_k )\T (d^{-1/2}w_j)\\
 &= \sqrt{\lambda_k(\Wc)}\sqrt{\lambda_k(\Wc) + \tau^2} v_{k j}(\Wc)   + O_p(d^{-1/4}).
\end{align*}
Since $d^{-1} w_j \T w_j = \sum_{\ell= 1}^m v^2_{\ell j}(\Wc) \lambda_\ell(\Wc)$, we have
$$r(\hat{w}_k,w_j) = \frac{d^{-1} \hat{w}_k\T w_j }{\| d^{-1/2}\hat{w}_k\|_2 \cdot\| d^{-1/2}w_j\|_2} \to v_{kj}(\Wc) { \zeta_{kj}^{1/2}}$$
in probability, as $d \to\infty$.

Proof of (ii). From Theorem~\ref{thm:main}, write
\begin{equation}\label{eq:proof3eq1}
d^{-1/2}\hat{w}_{k*} = \rho_k^{-1} \sum_{\ell = 1}^m v_{k\ell}(\Wc) d^{-1/2}  w_{\ell *} + O_p(d^{-1/2}),
 \end{equation}
 and note that $\E(w_{k*}) = \E(\hat{w}_{k*}) = 0 $. Then for $k =1,\ldots,m$, we have
$$
 \Var( d^{-1/2}  w_{k *}) = d^{-1} \E(w_{k *})^2 =  \sigma_{k}^2\E(z_{k *})^2 = \sigma_{k}^2,
$$
and, by (\ref{eq:proof3eq1}),
$$
 \Var( d^{-1/2}  \hat{w}_{k *} \mid W_1) = \rho_k^{-2} \sum_{\ell = 1}^m \left( v_{k\ell}(\Wc)\right)^2 \sigma_{\ell}^2 + O(d^{-1/2}).
$$
The independence of $w_{\ell *}$ and $w_{k *}$ for $k\neq \ell$ and (\ref{eq:proof3eq1}) give
\begin{align*}
\Cov(d^{-1/2}  \hat{w}_{k *}, d^{-1/2} w_{j *} \mid W_1)
 &=  \E( d^{-1}  \hat{w}_{k *}  w_{j *}  \mid W_1 )  \\
 &=  \rho_k^{-1}   v_{k j}(\Wc)  \sigma_j^2   + O(d^{-1/2}),
 \end{align*}
which in turn leads to
\begin{align*}
\mbox{corr}(\hat{w}_{k *},w_{j *} \mid W_1)
  &= \frac{\Cov(d^{-1/2}  \hat{w}_{k *}, d^{-1/2} w_{j *} \mid W_1) }{  \left( \Var( d^{-1/2}  w_{j *}) \Var( d^{-1/2}  \hat{w}_{k *} \mid W_1)  \right)^{1/2}} \\
  &= v_{k j}(\Wc)\frac{ {\sigma_j}}{\left[\sum_{\ell = 1}^m \left( v_{k\ell}(\Wc)\right)^2 \sigma_{\ell}^2\right]^{1/2} }   +  O(d^{-1/2}).
 \end{align*}
 \end{proof}


\section{Proof of Corollary 1}
\setcounter{equation}{0}
 \begin{corollary}\label{cor:rho_estimate} Suppose the assumptions of Lemma~\ref{thm:PCA-asymptotic} are satisfied. Let $d\to\infty$. For $i = 1,\ldots,m$, conditional to $W_1$, $\tilde\tau^2$, $\tilde\lambda_i(\Wc)$ and $\tilde{\rho}_i$ are consistent estimators of $\tau^2$, $\lambda_i(\Wc)$ and $\rho_i$, respectively.
\end{corollary}
\begin{proof}[Proof of Corollary~\ref{cor:rho_estimate}]
Lemma~\ref{thm:PCA-asymptotic} is used to show that $\tilde\tau^2$ and $\tilde\lambda_i(\Wc)$ converge in probability to $\tau^2$ and $\lambda_i(\Wc)$ as $d\to\infty$, respectively. By continuous mapping theorem, $\tilde\rho_i$ converges in probability to $\rho_i$.
\end{proof}

\section{Complete Table 2}

\begin{table}[t]
{\footnotesize
\centering
\begin{tabular}{cccccccc}
& & & & & $\rho_1$ \\
\cline{4-8}
   & $d$     & $n$   & Theory & Best & Asymp. & Jackknife & LZW \\
       \hline
 & 5000  & 50  & 1.41 (0.07) & 1.42 & 1.40 & 1.43 & 1.41 \\
Spike model & 10000 & 50  & 1.42 (0.06) & 1.43 & 1.42 & 1.44 & 1.42 \\
$\beta = 0.3$ & 10000 & 100 & 1.23 (0.03) & 1.23 & 1.23 & 1.24 & 1.23 \\
 & 20000 & 100 & 1.23 (0.02) & 1.23 & 1.23 & 1.24 & 1.23 \\
 \hline
 & 5000  & 50  & 1.42 (0.08) & 1.45 & 1.41 & 1.45 & 1.40 \\
Spike model & 10000 & 50  & 1.43 (0.07) & 1.45 & 1.43 & 1.46 & 1.42 \\
$\beta = 0.5$ & 10000 & 100 & 1.22 (0.02) & 1.23 & 1.22 & 1.23 & 1.21 \\
 & 20000 & 100 & 1.23 (0.02) & 1.23 & 1.23 & 1.24 & 1.22 \\
 \hline
& 5000  & 50  & 2.06 (0.06) & 2.22 & 1.92 & 2.14 & 2.00 \\
Mixture model & 10000 & 50  & 2.09 (0.06) & 2.17 & 1.98 & 2.14 & 2.02 \\
$a = 0.15$ & 10000 & 100 & 1.63 (0.02) & 1.67 & 1.61 & 1.65 & 1.63 \\
 & 20000 & 100 & 1.64 (0.02) & 1.66 & 1.62 & 1.66 & 1.63
\end{tabular}

\begin{tabular}{cccccccc}
& & & & & $\rho_2$ \\
\cline{4-8}
 & $d$     & $n$   & Theory & Best & Asymp. & Jackknife & LZW \\
 \hline
 & 5000  & 50  & 1.79 (0.11) & 1.86 & 1.75 & 1.78 & 1.79 \\
Spike model & 10000 & 50  & 1.79 (0.11) & 1.82 & 1.77 & 1.77 & 1.79 \\
$\beta = 0.3$ & 10000 & 100 & 1.43 (0.06) & 1.44 & 1.43 & 1.42 & 1.43 \\
 & 20000 & 100 & 1.43 (0.05) & 1.44 & 1.43 & 1.42 & 1.43 \\
 \hline
 & 5000  & 50  & 1.79 (0.11) & 1.99 & 1.72 & 1.81 & 1.71 \\
Spike model & 10000 & 50  & 1.80 (0.11) & 1.88 & 1.76 & 1.79 & 1.74 \\
$\beta = 0.5$ & 10000 & 100 & 1.44 (0.05) & 1.47 & 1.43 & 1.44 & 1.41 \\
 & 20000 & 100 & 1.42 (0.05) & 1.44 & 1.42 & 1.41 & 1.40 \\
 \hline
 & 5000  & 50  & 2.62 (0.21) & 5.44 & 2.20 & 2.68 & 2.46 \\
Mixture model & 10000 & 50  & 2.68 (0.19) & 3.20 & 2.35 & 2.68 & 2.50 \\
$a = 0.15$ & 10000 & 100 & 2.00 (0.09) & 2.13 & 1.90 & 2.00 & 1.99 \\
 & 20000 & 100 & 1.99 (0.10) & 2.05 & 1.93 & 1.97 & 1.97
\end{tabular}
}
\caption{Simulation results from 100 repetitions. ``Theory'' is mean (standard deviation) of $\rho_i$; ``Best'' is $\check\rho_i$ ; ``Asymp.'' is $\tilde\rho_i$ ; ``Jackknife'' is $\hat\rho_i^{(1)}$ ; ``LZW'' is from \cite{lee2010convergence}. Averages are shown for the latter four columns. The standard errors of the quantities in estimation of $\rho_i$ are at most 0.04.
  \label{tab:simulation}}
\end{table}

%

\bibliographystyle{biometrika}
\bibliography{library}

\end{document}